\documentstyle[twoside]{article}
\include{graphicx}
\title{\bf Refined asymptotics of the Riemann-Siegel theta function}
\author{\sc R. B. Paris\footnote{E-mail address:\ \ {\tt r.paris@abertay.ac.uk}}\\
\\
{\em Division of Computing and Mathematics,}\\
{\em Abertay University, Dundee DD1 1HG, UK}\\
}

\begin{document}
\newcommand{\bee}{\begin{equation}}
\newcommand{\ee}{\end{equation}}
\def\f#1#2{\mbox{${\textstyle \frac{#1}{#2}}$}}
\def\dfrac#1#2{\displaystyle{\frac{#1}{#2}}}
\newcommand{\fr}{\frac{1}{2}}
\newcommand{\fs}{\f{1}{2}}
\newcommand{\g}{\Gamma}
\newcommand{\la}{\lambda}
\newcommand{\al}{\alpha}
\newcommand{\br}{\biggr}
\newcommand{\bl}{\biggl}
\renewcommand{\topfraction}{0.9}
\renewcommand{\bottomfraction}{0.9}
\renewcommand{\textfraction}{0.05}
\newcommand{\gtwid}{\raisebox{-.8ex}{\mbox{$\stackrel{\textstyle >}{\sim}$}}}
\newcommand{\ltwid}{\raisebox{-.8ex}{\mbox{$\stackrel{\textstyle <}{\sim}$}}}
\newcommand{\mcol}{\multicolumn}
\date{}
\maketitle
\pagestyle{myheadings}
\markboth{\hfill {\it R.B. Paris} \hfill}
{\hfill {\it Riemann-Siegel theta function} \hfill}
\begin{abstract} 
The Riemann-Siegel theta function $\vartheta(t)$ is examined for $t\to+\infty$. Use of the refined asymptotic expansion for $\log\,\g(z)$ shows that the expansion of $\vartheta(t)$ contains an infinite sequence of increasingly subdominant exponential terms, each multiplied by an asymptotic series involving inverse powers of $\pi t$. Numerical examples are given to detect and confirm the presence of the first three of these exponentials. 

\vspace{0.4cm}

\vspace{0.3cm}

\end{abstract}

\vspace{0.2cm}

\noindent $\,$\hrulefill $\,$

\vspace{0.2cm}

\begin{center}
{\bf 1. \  Introduction}
\end{center}
\setcounter{section}{1}
\setcounter{equation}{0}
\renewcommand{\theequation}{\arabic{section}.\arabic{equation}}
The Riemann-Siegel theta function $\vartheta(t)$, which arises when the Riemann zeta function $\zeta(s)$ is expressed as the real-valued function $Z(t):=\exp (i\vartheta(t)) \zeta(\fs+it)$ on the critical line $\Re (s)=\fs$, is defined for real $t$ by 
\bee\label{e11}
\vartheta(t)=\arg\g(\fs it+\f{1}{4})-\fs t \log\,\pi.
\ee
This function satisfies $\vartheta(0)=0$ and is clearly an odd function of $t$, so that it is sufficient to consider $t>0$.
The well-known asymptotic expansion for $\vartheta(t)$ is given by \cite[p.~xiv]{Has}
\bee\label{e12}
\vartheta(t)\sim \frac{1}{2}t( \log\,(t/2\pi)-1)-\frac{\pi}{8}+\sum_{k=1}^\infty\frac{(1-2^{1-2k}) |B_{2k}|}{4k(2k-1) t^{2k-1}}+\frac{1}{2} \arctan (e^{-\pi t})
\ee
valid as $t\to+\infty$, where $B_{2k}$ denote the Bernoulli numbers. 

Recently, Brent \cite{RPB} employed an alternative representation of $\vartheta(t)$, obtained by application of the reflection and duplication formulas for the gamma function, in the form
\[\vartheta(t)=\frac{1}{2} \arg\,\g(\fs+it)-\frac{t}{2} \log\,2\pi-\frac{\pi}{8}+\frac{1}{2} \arctan (e^{-\pi t}).\]
By means of an improved error bound for the asymptotic expansion of $\log\,\g(z)$ valid in $\Re (z)\geq0$, $z\neq 0$, he was able to derive a rigorous error bound for the expansion in (\ref{e12}) valid for all $K\geq1$ and $t>0$ when the sum on the right-hand side of (\ref{e12}) is truncated after $K$ terms. 

The inclusion of the exponentially small term $\fs \arctan (e^{-\pi t})$ was shown to improve the attainable accuracy
of this bound. This term, which may be expressed for $t>0$ as a series of decreasing exponentials since
\bee\label{e12a}
\arctan (e^{-\pi t})=e^{-\pi t}-\f{1}{3}e^{-3\pi t}+\f{1}{5} e^{-5\pi t}-\f{1}{7} e^{-7\pi t}+ \cdots\ ,
\ee 
represents, as we shall show below, the first in a sequence of small exponentials present in the asymptotic expansion of $\vartheta(t)$.
A further application of the duplication formula enables $\vartheta(t)$ to be written in terms of $\g(it)$ and $\g(2it)$ as \cite{Gab}
\begin{eqnarray}
\vartheta(t)&=&\frac{1}{2}\Im\ \log\,\{\g(2it)/\g(it)\}-\frac{1}{2}t \log\,8\pi-\frac{\pi}{8}+\frac{1}{2} \arctan (e^{-\pi t})\nonumber\\
&=&\frac{1}{2}t(\log\,(t/2\pi)-1)+\frac{1}{2}\Im \{\Omega(2it)-\Omega(it)\}-\frac{\pi}{8}+\frac{1}{2} \arctan (e^{-\pi t}),\label{e13}
\end{eqnarray}
where
\bee\label{e14}
\log\,\g(it)=(it-\fs) \log\,it-it+\fs\log 2\pi+\Omega(it),
\ee
with $\Omega(it)$ denoting the slowly varying part of $\g(it)$.

The refined asymptotic expansion of $\log\,\g(z)$, first discussed in \cite{PW} and subsequently in \cite{B} and described in detail in \cite[\S 6.4]{PK}, can then be brought to bear on the asymptotic expansion of $\vartheta(t)$ for large $t>0$. It was shown \cite{B, PW} that the imaginary $z$-axis is a Stokes line for $\log\,\g(z)$, where in the neighbourhood of $\arg\,z=\pm\fs\pi$ an infinite number of increasingly subdominant exponential terms switch on as $|\arg\,z|$ increases. In this paper we shall exploit this theory to reveal the exponentially small contributions present in the expansion of $\vartheta(t)$. To achieve this, and to present a reasonably self-contained account, we reproduce the essential features of the refined asymptotics of $\log\,\g(z)$ when $\arg\,z=\fs\pi$ in Section 2. The details of our calculation for $\vartheta(t)$ are presented in Section 3 with a numerical verification given in Section 4.

\vspace{0.6cm}

\begin{center}
{\bf 2. \ The refined asymptotics of $\log\,\g(z)$ when $\arg\,z=\fs\pi$}
\end{center}
\setcounter{section}{2}
\setcounter{equation}{0}
\renewcommand{\theequation}{\arabic{section}.\arabic{equation}}
The slowly varying part of the  gamma function $\Omega(z)$ defined by $\g(z)=\sqrt{2\pi}z^{z-1/2}e^{-z+\Omega(z)}$ is shown in \cite[p.~282]{PK} to be given by the Mellin-Barnes integral
\[\Omega(z)=\frac{1}{2\pi i}\int_{c-\infty i}^{c+\infty i} \frac{\g(s)}{\sin \pi s}\,(2\pi z)^{-s} \zeta(1+s) \sin \fs\pi s\,ds\]
\[\hspace{0.6cm}=\sum_{k\geq 1}\frac{1}{k}\,\frac{1}{2\pi i}\int_{c-\infty i}^{c+\infty i} \frac{\g(s)}{\sin \pi s}\,(2\pi z)^{-s}  \sin \fs\pi s\,ds\]
where $|\arg\,z|\leq\pi-\delta$, $\delta>0$ and $0<c<1$. Displacement of the integration path to the right over the simple poles of the integrand at $s=1, 3, \ldots\,, 2n_k-1$, where $\{n_k\}_{k\geq1}$ is (for the moment) an arbitrary set of positive integers, and use of Cauchy's theorem then shows that, provided $|\arg\,z|\leq\pi-\delta$,
\bee\label{e20}
\Omega(z)=\frac{1}{\pi}\sum_{k\geq1}\frac{1}{k}\sum_{r=0}^{n_k-1} \frac{(-)^r (2r)!}{(2\pi kz)^{2r+1}}+{\cal R}(z;n_k)
%\sum_{k\geq1}\frac{R_k(z;n_k)}{k}.
\ee
With $L$ denoting the path parallel to the imaginary axis $\Re (s)=-c+2n_k+1$, $0<c<1$, the remainder term is given by (see \cite[p.~283]{PK})
\[{\cal R}(z;n_k)=\sum_{k\geq 1}\frac{1}{k}\,\frac{1}{2\pi i}\int_L \frac{\g(s)}{\sin \pi s}\,(2\pi kz)^{-s} \sin \fs\pi s\,ds\hspace{2cm}\]
\bee\label{e23a}
\hspace{2.8cm}=\sum_{k\geq1}\frac{1}{k}\bl\{e^{2\pi ikz}T_{2n_k+1}(2\pi ikz)-e^{-2\pi ikz} T_{2n_k+1}(-2\pi ikz)\br\},
\ee
where $T_\nu(z)$ is the so-called {\it terminant function} defined as a multiple of the incomplete gamma function $\g(a,z)$ by \cite[p.~242]{PK}
\bee\label{e21}
T_\nu(z)=\frac{e^{\pi\nu i}\g(\nu)}{2\pi i}\,\g(1-\nu,z)=\frac{e^{\pi i\nu} z^{-\nu}e^{-z}}{4\pi} \int_{-c-\infty i}^{-c+\infty i} \frac{\g(s+\nu)}{\sin \pi s}\,z^{-s}ds
\ee
with the integral being valid when $|\arg\,z|<\f{3}{2}\pi$ and $0<c<1$.
%the Mellin-Barnes integral representation \cite[p.~242]{PK}
%\[e^z T_\nu(z)=\frac{e^{\pi i\nu} z^{-\nu}}{4\pi} \int_{-c-\infty i}^{-c+\infty i} \frac{\g(s+\nu)}{\sin \pi s}\,z^{-s}ds\qquad(0<c<1,\ |\arg\,z|<\f{3}{2}\pi).\]
Thus we have the expansion
\bee\label{e23}
\Omega(z)=\frac{1}{\pi}\sum_{k\geq1}\frac{1}{k}\sum_{r=0}^{n_k-1} \frac{(-)^r (2r)!}{(2\pi kz)^{2r+1}}+
{\cal R}(z;n_k)
%\sum_{k\geq1}
%\frac{1}{k}\bl\{e^{2\pi ikz}T_{2n_k+1}(2\pi ikz)-e^{-2\pi ikz} T_{2n_k+1}(-2\pi ikz)\br\}
\ee
in the sector $|\arg\,z|\leq\pi-\delta$.
It is important to stress that the infinite sum on the right-hand side of (\ref{e23a}) is absolutely convergent (since
$\g(1-\nu,z)=O(z^{-\nu}e^{-z})$ for fixed $\nu$ and large $|z|$ in $|\arg\,z|\leq\f{3}{2}\pi-\delta$ \cite[(8.11.2)]{DLMF}) and that (\ref{e23}) is consequently {\it exact}.

Then, when $z=it$ ($t>0$), we obtain from (\ref{e23})
\bee\label{e22}
\Omega(it)=\frac{-i}{\pi}\sum_{k\geq 1}\frac{1}{k}\sum_{r=0}^{n_k-1}\frac{(2r)!}{(2\pi kt)^{2r+1}}+{\cal R}(it;n_k),
\ee
where
\bee\label{e22a}
{\cal R}(it;n_k)=\sum_{k\geq 1}\frac{1}{k}\bl\{e^{-2\pi kt} T_{2n_k+1}(2\pi kte^{\pi i})-e^{2\pi kt} T_{2n_k+1}(2\pi kt)\br\}
\ee
Two important features of the  expansion (\ref{e22}) are that: (i) the Stirling series appearing in the standard expansion of $\log\,\g(it)$ has been decomposed into a $k$-sequence of component asymptotic series with scale $2\pi kt$, each associated with its own arbitrary truncation index $n_k$, and (ii) the order of the terminant functions depends on $n_k$. It is these two features that permit each finite series in (\ref{e22}) to be optimally truncated  for large $t$ near its least term; that is, when
\bee\label{e24}
n_k=N_k\sim \pi kt.
\ee
When $n_k$ is chosen to be the optimal truncation index $N_k$ for the $k$th series, the order and argument of the  corresponding  terminant functions are approximately equal. To proceed further we now require the asymptotic expansions of $T_\nu(x)$ and $T_\nu(xe^{\pi i})$ when $\nu\sim x$ as $x\to+\infty$. 
\vspace{0.3cm}

\noindent{\bf 2.1\ \ Asymptotic expansions of $T_\nu(x)$ and $T_\nu(xe^{\pi i})$.}\ \  
\vspace{0.2cm}

\noindent The asymptotic expansion of $T_\nu(z)$ for large $|\nu|$ and $|z|$, when $|\nu|\sim|z|$, has been discussed in detail by Olver \cite{Olv}. His analysis was based on the Laplace integral representation
\bee\label{e25}
T_\nu(z)=\frac{e^{(\pi-\theta)\nu i}e^{-z}}{2\pi i}\int_0^\infty e^{-|z|t}\,\frac{t^{\nu-1}}{1+te^{-i\theta}}\,dt\qquad (z=xe^{i\theta},\ x>0,\ |\theta|<\pi).
\ee
Letting
\bee\label{e25a}
\nu=x+a,\qquad \nu>0,
\ee
where $a$ is bounded, we find that when $\theta=0$ the above integral is associated with a saddle point
at $t=1$ and may written as
\[T_\nu(x)=\frac{e^{\pi\nu i}e^{-2x}}{2\pi i} \int_{-\infty}^\infty e^{-\fr xw^2}\,\frac{t^{a-1}}{1+t}\,\frac{dt}{dw}\,dw,\qquad \fs w^2=t-\log\,t-1.\]
The branch of $w(t)$ is chosen such that $w\sim t-1$ as $t \to 1$ and 
reversion of the $w$--$t$ mapping yields
\[t=1+w+\f{1}{3}w^2+\f{1}{36}w^3-\f{1}{270}w^4+\f{1}{4320}w^5+ \cdots\ ,\]
from which we can compute, with the help of {\it Mathematica}, the series expansion
\[\frac{t^{a-1}}{1+t}\,\frac{dt}{dw}=\frac{1}{2}\sum_{s=0}^\infty {\bf A}_s(a) w^s.\]
Substitution of this expansion in the integral for $T_\nu(x)$ then yields \cite[\S 5]{Olv} 
\bee\label{e26}
T_\nu(x)\sim\frac{-ie^{\pi i\nu}e^{-2x}}{2\sqrt{2\pi x}}\sum_{s=0}^\infty {\bf A}_{2s}(a)\,\frac{\g(s+\fs)}{\g(\fs)}\, (\fs x)^{-s}\qquad(x\to+\infty),
\ee
where the first five even coefficients are
\begin{eqnarray}
{\bf A}_0(a)\!\!&=&\!\!1,\quad
{\bf A}_2(a)\!\!=\!\!\f{1}{6}(2-6a+3a^2),\nonumber\\
{\bf A}_4(a)\!\!&=&\!\!\f{1}{24cdot 6^2}(-11-120a+300a^2-192a^3+36a^4),\nonumber\\
{\bf A}_6(a)\!\!&=&\!\!\f{1}{300\cdot 6^4}(-587+3510a+9765a^2-26280a^3+18900a^4-5400a^5+540a^6),\label{e28a}\\
{\bf A}_8(a)\!\!&=&\!\!\f{1}{5600\cdot 6^6}(120341-44592a-521736a^2-722880a^3+2336040a^4-1826496a^5\nonumber\\
&&\hspace{5cm}+635040a^6-103680a^7+6480a^8),\nonumber\\
{\bf A}_{10}(a)\!\!&=&\!\!\f{1}{58800\cdot 6^8}(-4266772 - 14047026 a + 18368889 a^2 + 19144272 a^3 + 23661792 a^4 \nonumber\\
&&- 
 88817904 a^5 + 73929240 a^6 - 28921536 a^7 + 5987520 a^8 - 
 635040 a^9 + 27216 a^{10}).\nonumber
\end{eqnarray}

As $\theta\to\pi$, the integral (\ref{e25}) has a saddle almost coincident with the pole at $t=e^{i(\theta-\pi)}$
and may be continued analytically to the interval $\pi\leq\theta<2\pi$ provided the integration path is deformed to pass over the pole. When $\theta=\pi$, it is found that \cite[\S 5]{Olv}
\bee\label{e27}
T_{\nu}(xe^{\pi i})\sim\frac{1}{2}-\frac{i}{\sqrt{2\pi x}}\sum_{k=0}^\infty  {\bf B}_{2s}(a)\,\frac{\g(s+\fs)}{\g(\fs)}\, (\fs x)^{-s}.
\ee
The coefficients ${\bf B}_{s}(a)$ are computed from the expansion 
\[\frac{t^{a-1}}{1-t}\,\frac{dt}{dw}=-\frac{1}{w}+\sum_{k=0}^\infty {\bf B}_{s}(a)w^s,\]
where
it is found that the first five even-order coefficients are
%\footnote{There was a misprint in the first term in ${\hat G}_{6,j}$ in \cite{P13}, which appeared as $-3226$ instead of $-3626$. This was pointed out by T. Pudlik \cite{TP}. The correct value was used in the numerical calculations described in \cite{P13}.}
\begin{eqnarray}
{\bf B}_{0}(a)\!\!&=&\!\!\f{2}{3}-a,\qquad {\hat{\bf B}}_{2}=\f{1}{15\cdot 6^2}(46-225a+270a^2-90a^3), \nonumber\\
{\bf B}_{4}(a)\!\!&=&\!\!\f{1}{70\cdot 6^4}(230-3969a+11340a^2-11760a^3+5040a^4
-756a^5),\nonumber\\
{\bf B}_{6}(a)\!\!&=&\!\!\f{1}{350\cdot 6^6}(-3626-17781a+183330a^2-397530a^3+370440a^4
-170100a^5\nonumber\\
&&\hspace{7cm}+37800a^6-3240a^7),\nonumber\\
{\bf B}_{8}(a)\!\!&=&\!\!\f{1}{231000\cdot 6^8}(-4032746+43924815a+88280280a^2-743046480a^3\label{e28b}\\
&&+1353607200a^4-1160830440a^5+541870560a^6
-141134400a^7\nonumber\\
&&\hspace{6cm}+19245600a^8-1069200a^9),\nonumber\\
{\bf B}_{10}(a)\!\!&=&\!\!\f{1}{7007000\cdot 6^{10}}(502522570 + 1850358861 a - 12222960750 a^2 - 12894191310 a^3\nonumber\\
&& + 
   103403860560 a^4 - 167009778936 a^5 + 133973920080 a^6 - 
   62315613360 a^7 \nonumber\\
   &&+ 17552414880 a^8 - 2951348400 a^9 + 
   272432160 a^{10} - 10614240 a^{11}).\nonumber
\end{eqnarray}
\vspace{0.3cm}

\noindent{\bf 2.2\ \ The asymptotic expansion of $\Im\,\Omega(it)$.}\ \  
\vspace{0.2cm}

From (\ref{e22a}) and the expansions in (\ref{e26}) and (\ref{e27}), with the truncation indices $n_k$ chosen to be the optimal indices $N_k$ defined in (\ref{e24}), we obtain
\bee\label{e28}
\Im {\cal R}(it;N_k)\sim -\sum_{k\geq 1}\frac{e^{-2\pi kt}}{2\pi\sqrt{k^3t}}\sum_{s=0}^\infty {\bf C}_s(a_k) (\pi kt)^{-s},
\ee
where, from (\ref{e25a}) with $\nu=2N_k+1$ and $x=2\pi kt$,  
\bee\label{e210}
a_k:=2(N_k-\pi kt)+1. 
\ee
The coefficients ${\bf C}_s(a)$ are defined by 
\[{\bf C}_s(a):=\bl(\frac{1}{2}{\bf A}_{2s}(a)+{\bf B}_{2s}(a)\br)\,\frac{\g(s+\fs)}{\g(\fs)},\]
and from the values given in (\ref{e28a}) and (\ref{e28b}) we have
\begin{eqnarray}
{\bf C}_0(a)\!\!&=&\!\!\f{7}{6}-a,\quad {\bf C}_1(a)=\f{1}{30\cdot 6^2}(136-495a+405a^2-90a^3),\nonumber\\
{\bf C}_2(a)\!\!&=&\!\!\f{1}{1120\cdot 6^3}(-695-20538a+54180a^2-43680a^3+13860a^4-1512a^5),\nonumber\\
{\bf C}_3(a)\!\!&=&\!\!\f{1}{1120\cdot 6^5}(-15953 + 55929 a + 388395 a^2 - 949410 a^3 + 767340 a^4 - 
 283500 a^5 + 49140 a^6 - 3240 a^7),    \nonumber\\
{\bf C}_4(a)\!\!&=&\!\!\f{1}{70400\cdot 6^8}(170640893 + 21630510 a - 598217400 a^2 - 
  2559569760 a^3+ 6176233800 a^4  \nonumber\\
  &&- 5034007440 a^5 + 2026775520 a^6 - 
  436233600 a^7 + 48114000 a^8 - 2138400 a^9)\nonumber\\
{\bf C}_5(a)\!\!&=&\!\! \f{1}{25025\cdot 2^{11}\cdot 6^7}(-8649703370 - 28280511909 a + 27178306155 a^2 + 28170272130 a^3 \nonumber\\
&& + 
   154158404400 a^4   - 357524183016 a^5 + 292552139880 a^6 - 
   124352308080 a^7 \nonumber\\
&&+ 30395645280 a^8 - 4313509200 a^9 + 
   330810480 a^{10} - 10614240 a^{11}).\label{e210a}
\end{eqnarray}

Thus we finally obtain the desired expansion
\bee\label{e29}
\Im \Omega(it)\sim -\frac{1}{\pi}\sum_{k\geq 1}\frac{1}{k}\sum_{r=0}^{N_k-1} \frac{(2r)!}{(2\pi kt)^{2r+1}}-\frac{1}{2\pi t^{1/2}}\sum_{k\geq 1}\frac{e^{-2\pi kt}}{k^{3/2}} \sum_{s=0}^\infty {\bf C}_s(a_k) (\pi kt)^{-s}
\ee
as $t\to+\infty$, where $N_k\sim\pi kt$ is the optimal truncation index of the $k$th series in the first sum and the parameters $a_k$ appearing in the coefficients ${\bf C}_s(a_k)$ are determined by (\ref{e210}).

\vspace{0.6cm}

\begin{center}
{\bf 3. \ Asymptotic expansion of $\vartheta(t)$}
\end{center}
\setcounter{section}{3}
\setcounter{equation}{0}
\renewcommand{\theequation}{\arabic{section}.\arabic{equation}}
From (\ref{e13}) we require the large-$t$ expansion of the quantity
\[\Upsilon:=\frac{1}{2}\Im \{\Omega(2it)-\Omega(it)\}.\]
The series part of $\Upsilon$ is, from (\ref{e29}), given by
\bee\label{e30}
\frac{1}{2\pi}\sum_{k\geq1}\frac{1}{k}\bl\{\sum_{r=0}^{N_k-1}\frac{(2r)!}{(2\pi kt)^{2r+1}}-\sum_{r=0}^{2N_k-1}\frac{(2r)!}{(4\pi kt)^{2r+1}}\br\}.
\ee
The summation over $k$ may be carried out by a straightforward regrouping of the terms in the above absolutely convergent double series as described in \cite{B} and \cite[p.~286]{PK}:
\[\sum_{k\geq 1}\sum_{r=0}^{N_k-1}\frac{\alpha_r}{k^{2r+2}}=\sum_{r=0}^{N_1-1}\alpha_r\sum_{k\geq 1}\frac{1}{k^{2r+2}}
+\sum_{r=N_1}^{N_2-1}\alpha_r\sum_{k\geq 2}\frac{1}{k^{2r+2}}+\sum_{r=N_2}^{N_3-1}\alpha_r\sum_{k\geq 3}\frac{1}{k^{2r+2}}+\cdots\]
\[=\sum_{m=1}^\infty \sum_{r=N_{m-1}}^{N_m-1} \alpha_r \zeta(2r+2,m),\]
where $\alpha_r:=(2r)!/(2\pi t)^{2r+1}$ for the first series in (\ref{e30}), $N_0\equiv 0$ and $\zeta(s,m)=\sum_{k=0}^\infty (k+m)^{-s}$ is the Hurwitz zeta function.
This yields for the series contribution to $\Upsilon$ the result
\[\frac{1}{2\pi}\sum_{m=1}^\infty\bl\{\sum_{r=N_{m-1}}^{N_m-1}\frac{(2r)! \zeta(2r+2,m)}{(2\pi t)^{2r+1}}-\sum_{r=2N_{m-1}}^{2N_m-1} \frac{(2r)! \zeta(2r+2,m)}{(4\pi t)^{2r+1}}\br\},\]

The part of $\Upsilon$ involving terminant functions becomes, from (\ref{e29}),
\[\frac{1}{4\pi t^{1/2}}\bl\{\sum_{k\geq1}\frac{e^{-2\pi kt}}{k^{3/2}}\sum_{s=0}^\infty {\bf C}_s(a_k)(\pi kt)^{-s}-2^{-1/2}
\sum_{k\geq1}\frac{e^{-4\pi kt}}{k^{3/2}}\sum_{s=0}^\infty {\bf C}_s(a_{2k})(2\pi kt)^{-s}\br\}\]
\[=\frac{1}{4\pi t^{1/2}}\sum_{k\geq1}  \frac{(-)^{k-1}e^{-2\pi k t}}{k^{3/2}}\sum_{s=0}^\infty{\bf C}_s(a_k)(\pi k t)^{-s}.\]
Thus we obtain the following expansion theorem which is the main result of the paper:
\newtheorem{theorem}{Theorem}
\begin{theorem}$\!\!\!.$ The following expansion for the Riemann-Siegel theta function holds:
\[\vartheta(t)\sim \frac{1}{2}t(\log (t/2\pi)-1)-\frac{\pi}{8}+\frac{1}{2} \arctan(e^{-\pi t})\]
\[+\frac{1}{2\pi}\sum_{m=1}^\infty\bl\{\sum_{r=N_{m-1}}^{N_m-1}\frac{(2r)! \zeta(2r+2,m)}{(2\pi t)^{2r+1}}-\sum_{r=2N_{m-1}}^{2N_m-1} \frac{(2r)! \zeta(2r+2,m)}{(4\pi t)^{2r+1}}\br\}\]
\bee\label{e31}
+\frac{1}{4\pi t^{1/2}}\sum_{k\geq 1}  \frac{(-)^{k-1}e^{-2\pi k t}}{k^{3/2}}\sum_{s=0}^\infty{\bf C}_s(a_k)(\pi k t)^{-s}
\ee
as $t\to +\infty$, where $N_k\sim\pi kt$ is the optimal truncation index of the $k$th asymptotic series, with $N_0\equiv 0$. The first few coefficients ${\bf C}_s(a_k)$ are given in (\ref{e210a}) with the parameters $a_k$ specified in (\ref{e210}).
\end{theorem}

Thus we have established that there is an infinite number of exponentially small contributions to the asymptotics of $\vartheta(t)$ for large $t>0$ of the form $\exp\,[-2\pi kt]$, $k=1, 2, \ldots\,$, each multiplied by an asymptotic series in inverse powers of $\pi kt$. In addition, there is an infinite number of exponentials of the form $\exp\,[-(2k-1)\pi t]$ resulting from the expansion of $\arctan (e^{-\pi t})$ in (\ref{e13}).

We observe that when the series in $\Upsilon$ are both truncated at the index $n$ for all $k\geq 1$, the sum appearing in (\ref{e30}) becomes
\[\frac{1}{2\pi}\sum_{k\geq1}\frac{1}{k}\bl\{\sum_{r=0}^{n-1}\frac{(2r)!}{(2\pi kt)^{2r+1}}-\sum_{r=0}^{n-1}\frac{(2r)!}{(4\pi kt)^{2r+1}}\br\}=\frac{1}{2\pi}\sum_{k\geq1}\frac{1}{k}\sum_{r=0}^{n-1}\frac{(2r)!(1-2^{-2r-1})}{(2\pi kt)^{2r+1}}\]
\[=\frac{1}{2\pi}\sum_{r=1}^n\frac{(2r-2)!(1-2^{1-2r})}{(2\pi t)^{2r-1}}\,\zeta(2r)=\sum_{r=1}^n\frac{(1-2^{1-2r}) \,|B_{2r}|}{4r(2r-1) t^{2r-1}}\]
upon use of the identity connecting the zeta function of argument $2r$ to the Bernoulli numbers $B_{2r}$ given by $\zeta(2r)=(2\pi)^{2r} |B_{2r}|/(2 (2r)!)$. Thus, when the truncation index is chosen to be fixed and independent of $k$ the above sum reduces to the (truncated) standard series given in (\ref{e12}).

\vspace{0.6cm}

\begin{center}
{\bf 4. \ Numerical results}
\end{center}
\setcounter{section}{4}
\setcounter{equation}{0}
\renewcommand{\theequation}{\arabic{section}.\arabic{equation}}
In order to verify the presence of the exponentially small terms in our asymptotic formula for $\vartheta(t)$, we first subtract off the main expansion terms by defining
\[\Theta(t):=\vartheta(t)-\frac{1}{2}t(\log (t/2\pi)-1)+\frac{\pi}{8}-\frac{1}{2} \arctan(e^{-\pi t}).\]
The (truncated) asymptotic series and the subdominant exponential terms appearing in (\ref{e31}) are written as 
\[H_p(t):=\frac{1}{2\pi}\sum_{m=1}^p\bl\{\sum_{r=N_{m-1}}^{N_m-1}\frac{(2r)! \zeta(2r+2,m)}{(2\pi t)^{2r+1}}-\sum_{r=2N_{m-1}}^{2N_m-1} \frac{(2r)! \zeta(2r+2,m)}{(4\pi t)^{2r+1}}\br\},\]
\[E_k(t):=\frac{(-)^{k-1}e^{-2\pi kt}}{4\pi k(kt)^{1/2}}\sum_{s=0}^n {\bf C}_s(a_k)(\pi kt)^{-s},\]
where $p, k\geq 1$ and $n\geq0$ are integers. 

To detect the presence of the $K$th exponential $E_K(t)$ it is necessary to ``peel off'' from $\vartheta(t)$
all exponentials corresponding to $k<K$ and all larger terms in the asymptotic series $H_p(t)$.
In the case $p=1$, we define 
\[F_1(t):=\Theta(t)-H_1(t).\]
In Fig.~1(a) the terms of $H_1(t)$ are plotted against ordinal number $r$ when $t=5$; the jump arises because the first sum involves $N_1$ terms while the second sum involves $2N_1$ terms. It is seen that this sum contains terms that are much smaller than the first exponential since $e^{-10\pi}\simeq 2\times 10^{-14}$. Thus we expect that $F_1(t)\sim E_1(t)$ as a leading approximation.
In Table 1 we show the computed value of $F_1(t)$
compared with the value of the first subdominant exponential $E_1(t)$ for different truncation index $0\leq n\leq5$. 
It is seen that there is excellent agreement.
\begin{figure}[th]
	\begin{center}	{\tiny($a$)}\includegraphics[width=0.375\textwidth]{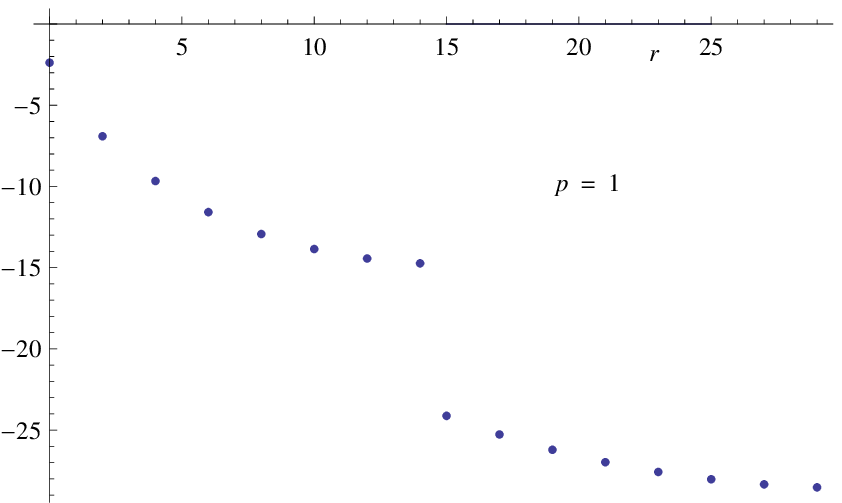}
	\qquad
	{\tiny($b$)}\includegraphics[width=0.375\textwidth]{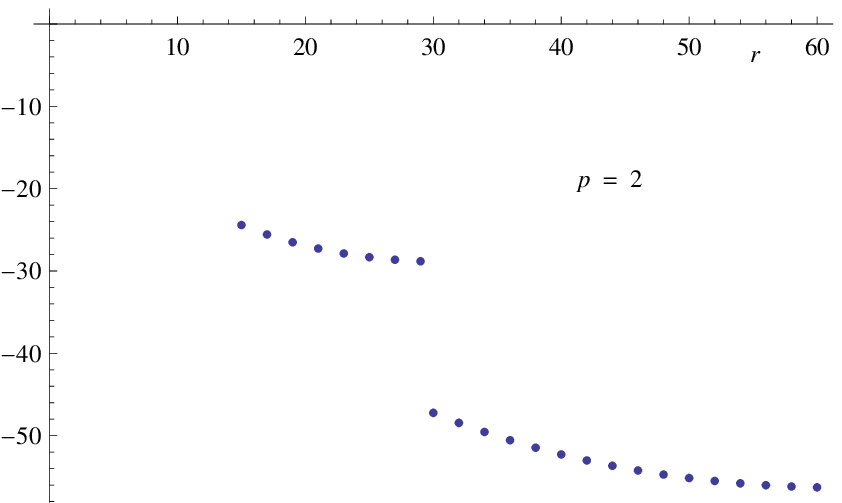}
\caption{\small{Magnitude of the terms (on a $\log_{10}$ scale) in $H_p(t)$ against ordinal number $r$ when $t=5$ and truncation indices $N_1=15$, $N_2=30$: (a) $p=1$ and (b) $p=2$.}}
	\end{center}
\end{figure}
\begin{table}[th]
\caption{\footnotesize{Values of the first subdominant exponential $E_1(t)$ as a function of truncation index $n$ for different $t$-values. The last line shows the value of $F_1(t)$.}} \label{t1}
\begin{center}
\begin{tabular}{|c|l|l|l|}
\hline
%&&\\[-0.3cm]
\mcol{1}{|c|}{} & \mcol{1}{c|}{$t=5,\ N_1=15$} & \mcol{1}{c|}{$t=8,\ N_1=25$} &  \mcol{1}{c|}{$t=10,\ N_1=31$} \\ 
\mcol{1}{|c|}{$n$} & \mcol{1}{c|}{$E_1(t)$} & \mcol{1}{c|}{$E_1(t)$} & \mcol{1}{c|}{$E_1(t)$} \\ 
\hline
&&&\\[-0.3cm]
0 & $1.2{\bf 7}911882757028\times 10^{-15}$ & $1.79{\bf 8}27918803825\times 10^{-24}$ & $1.29{\bf 6}04390329845\times 10^{-29}$ \\
1 & $1.299{\bf 0}5367126033\times 10^{-15}$ & $1.79141{\bf 9}07764499\times 10^{-24}$ & $1.2984{\bf8}420144033\times 10^{-29}$  \\
2 & $1.29933{\bf 5}26581188\times 10^{-15}$ & $1.791415{\bf 0}4299701\times 10^{-24}$ & $1.2984689{\bf 0}774446\times 10^{-29}$  \\
3 & $1.2993382{\bf 1}202001\times 10^{-15}$ & $1.791415{\bf 9}0147063\times 10^{-24}$ & $1.29846891{\bf 0}12541\times 10^{-29}$ \\
4 & $1.29933826{\bf 8}30296\times 10^{-15}$ & $1.791415887{\bf 6}1045\times 10^{-24}$ & $1.298468911{\bf8}1012\times 10^{-29}$ \\
5 & $1.299338269{\bf 2}8455\times 10^{-15}$ & $1.7914158875{\bf 7}452\times 10^{-24}$ & $1.298468911773{\bf4}0\times 10^{-29}$\\
[.1cm]\hline
&&&\\[-0.3cm]
$F_1(t)$ & $1.29933826977440\times 10^{-15}$ & $1.79141588758449\times 10^{-24}$ & $1.29846891177366\times 10^{-29}$\\
[.1cm]\hline
\end{tabular}
\end{center}
\end{table} 

To demonstrate the presence of the next two exponentials, we set
\[F_2(t):=\Theta(t)-H_2(t)-E_1(t)\]
and
\[F_3(t):=\Theta(t)-H_3(t)-E_1(t)+E_2(t).\]
In Fig.~1(b) the terms in $H_2(t)$ associated with $m=2$ are shown, from which it is seen that they decrease below the value $e^{-20\pi}\simeq 5\times 10^{-28}$; a similar conclusion arises (not shown) for the terms in $H_3(t)$ associated with $m=3$, where $e^{-30\pi}\sim1\times 10^{-41}$. Thus we expect the leading behaviours to satisfy $F_2(t)\sim -E_2(t)$ and $F_3(t)\sim E_3(t)$.

A difficulty arises in the computation of the exponential series appearing in $F_1(t)$ and $F_2(t)$, since these need
to be evaluated at optimal truncation (or at least as accurate as the following exponential series). With only the coefficients ${\bf C}_s(a)$ with $s\leq 5$ it is found that this is insufficient to achieve optimal truncation. To circumvent this problem, we computed the exponential terms in $F_2(t)$ and $F_3(t)$ from the terminant function representation in (\ref{e22a}) using the definition in terms of incomplete gamma functions in (\ref{e21}). The results are presented in Table 2, which clearly confirm the expansion in Theorem 1. 
\begin{table}[th]
\caption{\footnotesize{Values of the second and third subdominant exponentials as a function of truncation index $n$ when $t=5$ with truncation indices $N_1=15$, $N_2=30$, $N_3=45$. The last line shows the value of $F_j(t)$ ($j=2,3$).}} \label{t2}
\begin{center}
\begin{tabular}{|c|l|l|}
\hline
%&&\\[-0.3cm]
%\mcol{1}{|c|}{} & \mcol{1}{c|}{$t=5,\ N_1=15$} & \mcol{1}{c|}{$t=8,\ N_1=25$} &  \mcol{1}{c|}{$t=10,\ N_1=31$} \\ 
\mcol{1}{|c|}{$n$} & \mcol{1}{c|}{$E_2(t)$} & \mcol{1}{c|}{$E_3(t)$} \\ 
\hline
&&\\[-0.3cm]
0 & $-{\bf 1}.9459871830\times 10^{-29}$ & $3.{\bf 5}416763558\times 10^{-43}$ \\
1 & $-2.00{\bf 2}5102726\times 10^{-29}$ & $3.6{\bf 8}51106466\times 10^{-43}$ \\
2 & $-2.004{\bf 3}436323\times 10^{-29}$ & $3.691{\bf 4}163246\times 10^{-43}$ \\
3 & $-2.00440{\bf 0}3285\times 10^{-29}$ & $3.6916{\bf 7}53673\times 10^{-43}$ \\
4 & $-2.0044020{\bf 2}24\times 10^{-29}$ & $3.691685{\bf 4}260\times 10^{-43}$ \\
5 & $-2.00440207{\bf 2}3\times 10^{-29}$ & $3.6916858{\bf 0}15\times 10^{-43}$ \\
[.1cm]\hline
&&\\[-0.3cm]
$F_j(t)$ & $-2.0044020737\times 10^{-29}$ & $3.6916858157\times 10^{-43}$\\
[.1cm]\hline
\end{tabular}
\end{center}
\end{table} 

\end{document}